\documentclass[11pt,a4paper]{article}
\usepackage{times}

\title{La Baguette Math\'emagique}
\author{Julyan Cartwright}

\begin{document}

\maketitle

If you throw a needle or stick at random onto a floor ruled with parallel lines,  such as the cracks between floorboards or tiles, from the proportion of times that the stick lands crossing a crack you can estimate $\pi$; can we get $e$ as well? Yes, we can. 
All of these aspects have been discussed before, but I haven't seen them discussed in this way: that one can estimate both $\pi$ and $e$ with the same Buffon's needle experiment.

This algorithm for the computation of $\pi$ by Monte Carlo methods
is Buffon's needle, first described by the Count Buffon in 1733 in a report to the Acad\'emie des Sciences  \cite{buffon1733} and then discussed in a supplement to his great work of Natural History published in 1777 \cite{buffon1777}. 
For its implementation we need a floor 
with parallel lines; the cracks between floorboards, or floor tiles if we agree
to consider only the lines running in one direction.
The lines are a distance $a$ apart, onto which a pen
of length $l$ is thrown at random. 
The final position
of the pen on the floor can be given in terms of the distance $x$ of its
midpoint from the nearest line and the angle $\phi$ it makes with
the lines. For the pen to cross one of these lines, it must be
that
\begin{equation}
   x \leq l/2 \sin \phi
,\end{equation}
and the probability of this happening is
\begin{equation}
   p = \frac{1}{a\pi} \int_0^\pi l \sin \phi\, d\phi = \frac{2l}{a\pi}
.\end{equation}
This equation can be used to compute $\pi$ in the following way: we throw the 
pen at random onto the floor a large number of times $n$ and count the number 
of times $m$ it lands on a crack. Then $\pi$ is 
approximated by
\begin{equation}
   \pi \simeq \frac{2ln}{am}
.\end{equation}
As the French word for stick is baguette, some writers have mistakenly assumed that Buffon was advocating throwing bread
around, and this error is still extant. In fact, Buffon himself seems neither to have done any needle or stick throwing, nor to have considered  his thought experiment as an experimental approximation to $\pi$. Buffon was thinking of a general class of problems in geometric probability. 
``But if, instead of throwing into the air a round object such as a coin, one threw an object of another shape, such as a square Spanish pistole [a type of coin], or a needle, a stick, etc., the problem would demand a little more geometry, although in general it would always be possible to give its solution by comparisons of space, as we are going to demonstrate'', he wrote \cite{robertson1986}.
It seems to have been Laplace who first associated Buffon's idea with throwing a stick to obtain $\pi$ \cite{behrends2014}; Laplace also solved the slightly different problem using square tiles \cite{laplace1812,arnow1994} (see also \cite{schuster1974}). There are many more extensions: one can also consider convex bodies instead of needles \cite{aleman1997}; or flexible bodies, like noodles \cite{ramaley1969,waymire1994}.

As $n$ increases, better estimates of $\pi$ are obtained, albeit very slowly; since the standard error is $\epsilon=\sqrt{p(1-p)/n}$, we have $n\sim 1/\epsilon^2$; around a million throws give three significant figures.
Mario Lazzarini claimed in 1901 to have obtained $\pi$ to $6$ decimal places by
throwing a needle 3408 times \cite{lazzarini1901}, but as pointed out by Gridgeman \cite{gridgeman1960} and Badger \cite{badger1994}, he almost certainly cheated by
choosing the best moment to stop. We expect the estimate obtained by needle
throwing to cross the true value of $\pi$ infinitely many times as $n$
increases, so if we truncate near one of the crossings we get an extremely good
result, but we can only do this by knowing beforehand what the answer should be.
If one doesn't cheat, convergence is slow and Buffon's needle is not a quick way to estimate $\pi$.

On the other hand, it seems that nature may use the technique. Ants may use a Buffon needle type algorithm to
evaluate the size of potential nest sites \cite{mallon2000}.
They do this not by throwing anything, but by making use of the similar result 
that the estimated area of a plane is
inversely proportional to the number of intersections $N$ between two sets of lines, of total lengths $S$ and $L$,
randomly scattered on to it: thus $2SL/(\pi N)$. ``Scouts [i.e., ants] using such a Buffon's needle algorithm will assess nest area as inversely proportional to the number of intersections they make between a first set of pheromone-marked paths and a second set of census paths''  \cite{mallon2000}.

But we can obtain $e$ at the same time as we get $\pi$ from this experiment. We get $e$ from the result that if we 
sample random numbers uniformly from [0,1] until the sum first exceeds 1, the expected number of draws is  $e$.
Russell \cite{russell1991} indicates that the result comes from Gnedenko \cite{gnedenko}: ``an exercise in Gnedenko
...  asks the reader to show that, if $U_1$, $U_2$, $\ldots$ are iid uniformly on (0, 1), if $S_n = \sum_{i=1}^n U_i$  and if $N$  is the minimum value of $n$ for which $S_n > 1$, then $E(N) = e$.'' He adds, ``The proof is not difficult, and one approach is to solve the general problem of the expected number of draws for the sum to exceed $x$, with $0 \leq x \leq 1$. You get an integro-differential equation, the solution of which is $e^x$.''

So, to estimate both $\pi$ and $e$, throw the needle on floor, and at each throw, as well as recording whether or not it crosses a line, record its distance from the line nearer to you as a proportion of the distance to the next line --- any point of the needle will do, an end or the middle, as long as it is consistent. 
Then $\pi$ is given by the proportion of times the needle crosses a line. And $e$ is given by the average number of observations that need to be added together to reach one.

\bibliographystyle{unsrt}
\bibliography{buffon}

\end{document}